\documentclass[11pt,a4paper,oneside]{article}

\usepackage[T1]{fontenc}
\usepackage[latin1]{inputenc}

\usepackage{amsmath}
\usepackage{amsfonts}
\usepackage{amssymb}
\usepackage{amsthm}
\usepackage{stmaryrd}
\usepackage{graphicx}
\usepackage{floatflt}

\usepackage[all,cmtip,color,arc,curve]{xy}

\title{On the Birational Geometry of Schubert Varieties}
\author{Benjamin Schmidt}

\theoremstyle{definition}
\newtheorem{defi}{Definition}[section]
\newtheorem{ex}[defi]{Example}

\theoremstyle{plain}
\newtheorem{thm}[defi]{Theorem}
\newtheorem{lem}[defi]{Lemma}
\newtheorem{prop}[defi]{Proposition}
\newtheorem{cor}[defi]{Corollary}

\begin{document}

\maketitle

\begin{abstract}
We classify all $\mathbb{Q}$-factorializations of (co)minuscule Schubert
varieties by using their Mori dream space structure. As a corollary we obtain a
description of all IH-small resolutions of (co)minuscule Schubert varieties
generalizing results of Perrin. We improve his results by including
algebraically closed fields of positive characteristic and cominuscule Schubert
varieties. Moreover, the use of $\mathbb{Q}$-factorializations and Mori dream
spaces simplifies the arguments substantially.
\end{abstract}

\section{Introduction}

A fundamental goal of algebraic geometry is to describe birational models with
better properties than the average variety. These models should be smooth or at least have mild singularities. A major step towards this goal was the resolution of singularities for any field of characteristic $0$ by Hironaka in \cite{Hir64}. For Schubert varieties there are the well known Bott-Samelson resolutions introduced in \cite{BS58}. They are rational resolutions and lead to a character formula for representations of reductive groups. On the other hand, Zelevinsky constructed IH-small resolutions (see Definition \ref{IHsmall}) in the case of Schubert varieties in Grassmannians and used them to compute Kazhdan-Lusztig polynomials (see \cite{Zel83}). Sankaran and Vanchinathan obtained similar results in Lagrangian and maximal isotropic Grassmannians in \cite{SV94} and \cite{SV95}.

Many classical results on Schubert varieties in Grassmannians can be generalized to minuscule or cominuscule Schubert varieties (see Definition \ref{defMinuscule}). In \cite{Per07} Perrin gives a complete classification of all IH-small resolutions of minuscule Schubert varieties over $\mathbb{C}$. This was done using a connection to the minimal model program: Totaro proved that any IH-small resolution is a relative minimal model in \cite[Proposition 8.3]{Tot00}. Perrin was able to classify all relative minimal models of minuscule Schubert varieties. Since most of the results from the minimal model program are only known in characteristic $0$, this approach is only valid over the complex numbers.

We will investigate further into the birational geometry of Schubert varieties
over an arbitrary algebraically closed field using a different approach. A
Mori-small morphism from a normal and $\mathbb{Q}$-factorial variety to a normal
variety is called a $\mathbb{Q}$-factorialization. Our goal is to determine all
$\mathbb{Q}$-factorializations of any (co)minuscule Schubert varieties.

In order to handle the occurring combinatorics in the Weyl group one defines a
quiver for each reduced expression (see Definition \ref{defQuiver}). Due to a
result in \cite{Ste96}, every reduced expression of a (co)minuscule element is
unique up to commuting relations. This implies that there is a unique quiver
associated to each (co)minuscule Schubert variety. Moreover, there is an
explicit combinatorial description of the quivers of minuscule elements (see
Theorem \ref{propQuiv2}). This provides a very concrete object to work with. We
define a partial ordering on the vertices of the quiver, call the maximal
elements peaks and assign each vertex a value called the height.

For each ordering of the peaks, there is a birational projective morphism
$\widehat{\pi}: \widehat{X}(\widehat{w}) \to X(w)$ (see Section \ref{prelim}).
The varieties $\widehat{X}(\widehat{w})$ are towers of locally trivial
fibrations with fibers being Schubert varieties. This generalizes the
Bott-Samelson resolution which is a tower of locally trivial
$\mathbb{P}^1$-fibrations. These varieties $\widehat{X}(\widehat{w})$ are
generally not smooth, but in our case always locally $\mathbb{Q}$-factorial and
normal. We prove the following theorem.

\newtheorem*{thm:qfact}{Theorem \ref{qfact}}
\begin{thm:qfact}
Let $X(w)$ be a (co)minuscule Schubert variety. Then all
$\mathbb{Q}$-factorializations of $X(w)$ are given by the morphisms
$\widehat{\pi}: \widehat{X}(\widehat{w}) \to X(w)$ obtained from any ordering of
the peaks.
\end{thm:qfact}

The use of Mori dream spaces is the main ingredient for proving this result.
These spaces are tailor-made for running the minimal model program in any
characteristic (see \cite{HK00}). First, we show that the varieties
$\widehat{X}(\widehat{w})$ are indeed Mori dream spaces. More precisely, taking
all $\widehat{X}(\widehat{w})$ for any possible ordering of the peaks leads to
all the small $\mathbb{Q}$-factorial modifications defining a Mori dream space.
Following \cite{Dem74} and \cite{Per07}, we give explicit descriptions of the
nef and effective cones of divisors using the structure of towers of locally
trivial fibrations. This is all that is needed to describe the Mori dream space
structure. The Theorem follows because the morphisms $\widehat{\pi}:
\widehat{X}(\widehat{w}) \to X(w)$ are all Mori-small, i.e., they do not
contract any divisor.

Since any IH-small morphism is also Mori-small, classifying all IH-small
resolutions of (co)minuscule Schubert varieties over any algebraically closed
field becomes a matter of checking which of these $\mathbb{Q}$-factorializations
are IH-small. Zelevinsky for Grassmannians and Perrin for minuscule homogeneous
spaces define specific orderings using heights of peaks that are called neat.
Generalizing their results by including cominuscule Schubert varieties and
algebraically closed fields of positive characteristic, we obtain the following
corollary.

\newtheorem*{cor:mainThm}{Corollary \ref{mainThm}}
\begin{cor:mainThm}
Let $X(w)$ be a (co)minuscule Schubert variety over an algebraically closed
field $k$. Then the IH-small resolutions of $X(w)$ are exactly given by the
morphisms $\widehat{\pi}: \widehat{X}(\widehat{w}) \to X(w)$, where
$\widehat{w}$ is obtained from a neat ordering of the peaks and
$\widehat{X}(\widehat{w})$ is smooth.
\end{cor:mainThm}

Note that there is an explicit combinatorial criterion for smoothness of the
varieties $\widehat{X}(\widehat{w})$ (see Section \ref{IHresolution}). Using
Theorem \ref{qfact} allows for a uniform treatment of both the minuscule and
cominuscule case. While similar arguments as in \cite{Per07} can be used at
least over $\mathbb{C}$, the geometry is slightly different. This results in
even more complicated combinatorics (see \cite{Sch11}). Therefore, we strongly
believe that the present proof is much more suitable to approach non
(co)minuscule cases.

Section 2 sets up some basic notation. In Section 3 we recall (co)minuscule
Schubert varieties and the definition of $\widehat{X}(\widehat{w})$. Section 4
is concerned with the proof of Theorem \ref{qfact}, while Section 5 presents
Corollary \ref{mainThm}. 

\paragraph*{Acknowledgments.} I would like to thank Nicolas Perrin and Emanuele
Macri for carefully reading preliminary versions of this manuscript. I am also
thanking Nicolas Perrin for many hours of explanation and Ana-Maria Castravet
for giving me useful hints on Mori dream spaces. I also appreciate useful
suggestions from the referee. This research was partially supported by the Hausdorff Center
for Mathematics in Bonn.

\section{Notation}

Let $G$ be a simple algebraic group over an algebraically closed field $k$. By
$T$ we denote a maximal torus in $G$ and $B$ is a Borel subgroup containing $T$.
The variety $X$ is usually the homogeneous space $G/P$ for a maximal parabolic
subgroup $P$. Furthermore, $W$ shall be the Weyl group of $G$ and $R$ the set of
all roots, while $S$ is the set of simple roots corresponding to $(B,T)$. We
denote the set of positive roots by $R^{+}$, while $R^{-}$ is the set of
negative roots. Moreover, $l$ is the length function on $W$ corresponding to
$S$. For root systems we use the notation from \cite{Bou68}.

For any projective normal $\mathbb{Q}$-factorial variety $Y$, we denote the cone
of effective $\mathbb{Q}$-divisors by $\operatorname{Eff}(Y)$ and the closed
cone of nef $\mathbb{Q}$-divisors by $\operatorname{Nef}(Y)$. The movable cone
$\operatorname{Mov}(Y)$ is the cone generated by the divisors with stable
base locus of codimension bigger than $1$. The inclusions
$\operatorname{Nef}(Y) \subset \overline{\operatorname{Mov}}(Y)$
and $\overline{\operatorname{Mov}}(Y) \subset
\overline{\operatorname{Eff}}(Y)$ follow directly from the definitions.

\section{Preliminaries}
\label{prelim}

In this section we are going to recall quivers corresponding to reduced expressions and define varieties $\widehat{X}(\widehat{w})$ from a decomposition of these quivers.

\subsection{Quivers}

Let $\widetilde{w} = (s_{\beta_1}, \ldots, s_{\beta_r})$ be a reduced decomposition of an element $w \in W$, i.e., $r$ is minimal such that $w = s_{\beta_1} \cdots s_{\beta_r}$ where $s_{\beta_1}, \ldots, s_{\beta_r}$ are simple reflections corresponding to $\beta_1, \ldots, \beta_r \in S$. Whenever it exists, the successor $s(i)$ of an index $i \in [1,r]$ is the smallest index $j>i$ such that $\beta_i = \beta_j$. Similarly, the predecessor $p(i)$ of an index $i \in [1,r]$ is the biggest index $j<i$ such that $\beta_i = \beta_j$. The combinatorics in the Weyl group can be translated into the geometry of the following quiver.

\begin{defi}
\label{defQuiver}
The quiver $Q_{\widetilde{w}}$ has vertices given by $[1,r]$. There is an arrow from $i$ to $j$ if $\langle \beta_i^{\vee}, \beta_j \rangle \neq 0$ and $i<j<s(i)$ (or $i<j$ if $s(i)$ does not exists). In addition, each vertex has a color via the map $[1,r] \to S$ given by $i \mapsto \beta_i$. A partial ordering on $Q_{\widetilde{w}}$ is generated by the relations $i \succeq j$ whenever there is an arrow from $i$ to $j$.
\end{defi}

Note, that the partial ordering is not the one defined in \cite{Per07}, but the reversed one. This is more natural with respect to all our notation and pictures. The quiver describes a reduced expression up to commuting relations between the simple reflections.

\begin{defi}
\label{defMinuscule}
Let $\omega$ be a fundamental weight corresponding to a simple root $\alpha$.
\begin{enumerate}
\item We call $\omega$ \textbf{minuscule} if $\langle \beta^{\vee}, \omega \rangle \leq 1$ for all $\beta \in R^+$.
\item We call $\omega$ \textbf{cominuscule} if the fundamental weight $\omega^{\vee}$ corresponding to $\alpha^{\vee} \in R^{\vee}$ is minuscule.
\end{enumerate}
\end{defi}

Equivalently, a fundamental weight $\omega$ corresponding to a simple root $\alpha$ is cominuscule if the coefficient of $\alpha$ at the highest root of $R$ is $1$. This characterization leads easily to a complete table of all the minuscule and cominuscule weights. \\
\\
\centerline{
\begin{tabular}[ht]{c|c|c}
Type & Minuscule weights & Cominuscule weights \\
\hline
$A_n$ & $\omega_1$, \ldots, $\omega_n$ & $\omega_1$, \ldots, $\omega_n$ \\
$B_n$ & $\omega_n$ & $\omega_1$ \\
$C_n$ & $\omega_1$ & $\omega_n$ \\
$D_n$ & $\omega_1$, $\omega_{n-1}$, $\omega_n$ & $\omega_1$, $\omega_{n-1}$, $\omega_n$ \\
$E_6$ & $\omega_1$, $\omega_6$ & $\omega_1$, $\omega_6$ \\
$E_7$ & $\omega_7$ & $\omega_7$ \\
$E_8$ & none & none \\
$F_4$ & none & none \\
$G_2$ & none & none \\
\end{tabular}}

For any (co)minuscule fundamental weight $\omega$ we define $P_{\omega}$ as the
maximal parabolic subgroup corresponding to $\omega$. An element $w \in W$ is
called (co)minuscule with respect to $\omega$ if it is a representative
of minimal length in its class modulo the Weyl group of $P_{\omega}$ denoted by
$W_{P_{\omega}}$. In this case the Schubert variety
$X_{P_{\omega}}(w):=\overline{BwP_{\omega}/P_{\omega}}$ is also called (co)minuscule. By \cite{Ste96} every (co)minuscule element has a unique reduced expression up to commuting relations. Therefore, we will usually write $Q_w$ instead of $Q_{\widetilde{w}}$. The next  theorem describes the shape of minuscule quivers. It is proven in \cite[Proposition 4.1]{Per07}. The fact that we do not exclude the non-simply laced case hardly changes anything in the proof. As a semisimple linear algebraic group and its dual group have the same Weyl group, the cominuscule case follows immediately.

\begin{thm}
\label{propQuiv2}
Let $\widetilde{w} = (s_{\beta_1}, \ldots, s_{\beta_r})$ be a reduced expression of an element $w$ in the Weyl group $W$ and $\omega$ a fundamental weight. Then $w$ is minuscule with respect to $\omega$ if and only if the following three conditions hold.
\begin{enumerate}
\item The element $\beta_r$ is the unique simple root such that $\langle \beta_r^{\vee}, \omega \rangle = 1$.
\item Let $i \succ r$ be a vertex of the quiver such that $s(i)$ does not exist. Then there is a unique arrow from $i$ to a vertex $k$ and we have $\langle \beta_i^{\vee}, \beta_k \rangle = -1$.
\item Let $i \succ r$ be a vertex of the quiver such that $s(i)$ exists. Then there are two possibilities. Either there are two distinct vertices $k_1$, $k_2$ with an arrow coming from $i$ or there is a unique vertex $k$ with an arrow coming from $i$. In the first case, 
$$\langle \beta_i^{\vee}, \beta_{k_1} \rangle = \langle \beta_i^{\vee}, \beta_{k_2} \rangle = -1.$$
In the second case,
$$\langle \beta_i^{\vee}, \beta_k \rangle = -2.$$
\end{enumerate}
\end{thm}

In the (co)minuscule case the Bruhat order can easily be described on the quiver. Let $w \in W$ be any (co)minuscule element. Then we get a subquiver $Q_{w'}$ of $Q_w$ by removing a maximal vertex. This subquiver corresponds to an element $w' \in W$ and $w' \leq w$ in the Bruhat order. In fact, these relations generate the Bruhat order (see \cite[Theorem 3.6]{Per07}). This means the Bruhat order equals the weak Bruhat order. Therefore, in order to understand the quivers it suffices to understand the quivers of the maximal elements. This can be obtained via the last theorem. They are simply given by the maximal quivers satisfying the three conditions.

Next, we want to decompose the quiver and define a variety that the Bott-Samelson resolution factorizes through. Let $w \in W$ be a (co)minuscule element. A peak of $Q_w$ is a maximal vertex and the set of peaks is denoted by $p(Q_w)$.

\begin{defi}
\begin{enumerate}
\item Let $p$ be a peak in $Q_w$. Then we define the set
$$\widehat{Q}_w(p) := \{ i \in Q_w \ | \ \exists q \in p(Q_w)\backslash \{p\}, \ i \preceq q \}.$$
The complement is denoted by $Q_w (p) := Q_w \backslash \widehat{Q}_w(p)$.
\item Let $p_1, \ldots p_s$ be any ordering of the peaks of $Q_w$. Inductively we define $Q_0 := Q_w$ and $Q_i := \widehat{Q}_{i-1}(p_i)$ for all $i \in [1,s]$.
\end{enumerate}
\end{defi}

By setting $Q_{w_i}:= Q_{i-1}(p_i)$ we obtain a generalized reduced
decomposition $\widehat{w}=(w_1, \ldots, w_s)$ of $w$, where $w_i$ is the
element of the Weyl group corresponding to the quiver $Q_{w_i}$. This means
$w=w_1 \cdots w_s$ and $l(w) = \sum_{i=1}^s l(w_i)$. The decomposition is made
in a way such that every $Q_{w_i}$ has a unique peak. Moreover, we write
$m_{\widehat{w}}(Q_w)$ for the union of all the minimal vertices of $Q_{w_i}$
for $i \in [1,s]$

\begin{ex}
Let us introduce a few examples. Later, we will further investigate them. The
following three diagrams show the decompositions $w_1 = (s_{\alpha_3}
s_{\alpha_4}) (s_{\alpha_1}s_{\alpha_2}s_{\alpha_3}s_{\alpha_4})$ in type $C_4$,
$w_2 = s_{\alpha_3} (s_{\alpha_1}s_{\alpha_2})
(s_{\alpha_5}s_{\alpha_4}s_{\alpha_3})$ in type $A_5$ and $w_3 =
(s_{\alpha_5}s_{\alpha_4}s_{\alpha_2})
(s_{\alpha_1}s_{\alpha_3}s_{\alpha_4}s_{\alpha_5}s_{\alpha_6})$ in type $E_6$.
We do not draw the direction of the arrows because they all go down. The
coloring is given by the projection onto the $C_4$, $A_5$, and $E_6$ quivers.

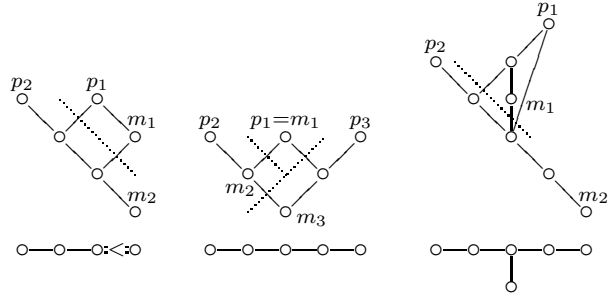
\centerline{
\xygraph{
!{<0cm,0cm>;<1cm,0cm>:<0cm,1cm>::}
!{(0,-1) }*{\circ}="1"
!{(1,-1) }*{\circ}="2"
!{(0.5,-1.5) }*{\circ}="3"
!{(1.5,-1.5) }*{\circ}="4"
!{(1,-2) }*{\circ}="5"
!{(1.5,-2.5) }*{\circ}="6"
!{(0,-3) }*{\circ}="a"
!{(0.5,-3) }*{\circ}="b"
!{(1,-3) }*{\circ}="c"
!{(1.5,-3) }*{\circ}="d"
!{(0.5,-1) }*{}="cut1"
!{(1.5,-2) }*{}="cut2"
!{(1,-0.8) }*{\scriptstyle{p_1}}="p1"
!{(0,-0.8) }*{\scriptstyle{p_2}}="p2"
!{(1.6,-1.3) }*{\scriptstyle{m_1}}="m1"
!{(1.6,-2.3) }*{\scriptstyle{m_2}}="m2"
!{(2.5,-1.5) }*{\circ}="1'"
!{(3.5,-1.5) }*{\circ}="2'"
!{(4.5,-1.5) }*{\circ}="3'"
!{(3,-2) }*{\circ}="4'"
!{(4,-2) }*{\circ}="5'"
!{(3.5,-2.5) }*{\circ}="6'"
!{(2.5,-3) }*{\circ}="a'"
!{(3,-3) }*{\circ}="b'"
!{(3.5,-3) }*{\circ}="c'"
!{(4,-3) }*{\circ}="d'"
!{(4.5,-3) }*{\circ}="e'"
!{(3,-1.5) }*{}="cut1'"
!{(3.5,-2) }*{}="cut2'"
!{(4,-1.5) }*{}="cut3'"
!{(3,-2.5) }*{}="cut4'"
!{(3.5,-1.3) }*{\scriptstyle{p_1=m_1}}="p1'"
!{(2.5,-1.3) }*{\scriptstyle{p_2}}="p2'"
!{(4.5,-1.3) }*{\scriptstyle{p_3}}="p3'"
!{(2.9,-2.2) }*{\scriptstyle{m_2}}="m2'"
!{(3.85,-2.6) }*{\scriptstyle{m_3}}="m3'"
!{(7,0) }*{\circ}="1''"
!{(5.5,-0.5) }*{\circ}="2''"
!{(6.5,-0.5) }*{\circ}="3''"
!{(6,-1) }*{\circ}="4''"
!{(6.5,-1) }*{\circ}="5''"
!{(6.5,-1.5) }*{\circ}="6''"
!{(7,-2) }*{\circ}="7''"
!{(7.5,-2.5) }*{\circ}="8''"
!{(5.5,-3) }*{\circ}="a''"
!{(6,-3) }*{\circ}="b''"
!{(6.5,-3) }*{\circ}="c''"
!{(7,-3) }*{\circ}="d''"
!{(7.5,-3) }*{\circ}="e''"
!{(6.5,-3.5) }*{\circ}="f''"
!{(5.75,-0.5) }*{}="cut1''"
!{(6.75,-1.5) }*{}="cut2''"
!{(7,0.2) }*{\scriptstyle{p_1}}="p1''"
!{(5.5,-0.3) }*{\scriptstyle{p_2}}="p2''"
!{(6.9,-1.1) }*{\scriptstyle{m_1}}="m1''"
!{(7.6,-2.3) }*{\scriptstyle{m_2}}="m2''"
"cut1"-@{.}"cut2"
"a"-"b"
"b"-"c"
"c"-@{=}|{<}"d"
"1"-"3"
"2"-"3"
"2"-"4"
"3"-"5"
"4"-"5"
"5"-"6"
"cut1'"-@{.}"cut2'"
"cut3'"-@{.}"cut4'"
"a'"-"b'"
"b'"-"c'"
"c'"-"d'"
"d'"-"e'"
"1'"-"4'"
"2'"-"4'"
"2'"-"5'"
"3'"-"5'"
"4'"-"6'"
"5'"-"6'"
"cut1''"-@{.}"cut2''"
"a''"-"b''"
"b''"-"c''"
"c''"-"d''"
"c''"-"f''"
"d''"-"e''"
"1''"-"3''"
"1''"-"6''"
"2''"-"4''"
"3''"-"4''"
"3''"-"5''"
"4''"-"6''"
"5''"-"6''"
"6''"-"7''"
"7''"-"8''"
}}
\end{ex}

\subsection{Intermediate Varieties}

We will now construct the \textbf{intermediate variety}
$\widehat{X}(\widehat{w})$ corresponding to one of the generalized reduced
decomposition $\widehat{w}$ described before. Let $P_{\beta_i}$ be the minimal
parabolic subgroup $B \cup Bs_{\beta_i}B$ corresponding to $\beta_i$ for $i \in
[1,r]$ and $p_{\beta_i}: G/B \to G/P_{\beta_i}$ be the quotient morphism. Its
fibers are isomorphic to $P_{\beta_i}/B \cong \mathbb{P}^1$. For any $x \in G/B$
we define $\mathbb{P}(x, \beta_i)$ to be $p_{\beta_i}^{-1}(p_{\beta_i}(x))$. Let
$P^{\beta_i}$ be the maximal parabolic subgroup containing $B$ corresponding to
$\beta_i$ for $i \in [1,r]$. For any $x \in G/B$ the map $G/B \to G/P^{\beta_i}$
restricted to $\mathbb{P}(x,\beta_i)$ is an isomorphism onto its image which is
called $\overline{\mathbb{P}}(x,\beta_i)$. If $y \in G/B$ and $\overline{x} \in
\overline{\mathbb{P}}(y,\beta_{i-1})$ for $x \in \mathbb{P}(y,\beta_{i-1})$ we
abuse notation by writing $\overline{\mathbb{P}}(\overline{x},\beta_i) =
\overline{\mathbb{P}}(x,\beta_i)$.

\begin{defi}
\begin{enumerate}
\item The Bott-Samelson variety $\widetilde{X}(\widetilde{w})$ can be defined as
(see \cite[Remark 2.7]{Per07})
$$\left\{ (x_1, \ldots, x_r) \in \prod_{i=1}^r G/P^{\beta_i} \ | \ x_0 = 1, \
x_i \in \overline{\mathbb{P}}(x_{i-1},\beta_i) \text{ for all } i \in [1,r]
\right\}.$$
\item We define $m_{\widehat{w}}(Q_w)$ to be the union of all the minimal
vertices of $Q_{w_i}$ for $i \in [1,s]$.
\item The variety $\widehat{X}(\widehat{w})$ is defined to be the image of
$\widetilde{X}(\widetilde{w})$ under the projection $\prod_{i=1}^r G/P^{\beta_i}
\to \prod_{i \in m_{\widehat{w}}(Q_w)} G/P^{\beta_i}$.
\end{enumerate}
\end{defi}

As the Bott-Samelson resolution is given by the projection to the last factor, it is clear that it factorizes through the intermediate variety.\\
\\
\centerline{
\xymatrix{
\widetilde{X}(\widetilde{w}) \ar[r]^{\widetilde{\pi}} \ar[rd]^{\pi} & \widehat{X}(\widehat{w}) \ar[d]^{\widehat{\pi}} \\
& X(w)
}} \\
\\

These intermediate varieties have some very nice properties as proven in \cite[Section 5]{Per07}. For the special decompositions of the quiver described above it is possible to define them as towers of locally trivial fibrations with fibers being the Schubert varieties corresponding to $w_i$ as follows. Let $m_i$ be the minimal vertex of the quiver $Q_{w_i}$ for $i \in [1,s]$. There are projection $f_i: \prod_{j=1}^i G/P^{\beta_{m_j}} \to \prod_{j=1}^{i-1} G/P^{\beta_{m_j}}$. 

\begin{thm} The restriction of $f_i$ to the image $f_{i+1} \circ \ldots \circ f_{s} (\widehat{X}(\widehat{w}))$ is a locally trivial fibration and the fiber is given by the Schubert variety corresponding to $w_i$. In particular, $\widehat{X}(\widehat{w})$ is normal.
\end{thm}

We can also understand divisors and the Picard group of the intermediate varieties very well. In order to explain this, we need to go back to the Picard group of the Bott-Samelson variety. For any $i \in [1,r]$ we define a divisor of the Bott-Samelson variety $\widetilde{X}(\widetilde{w})$ by
$$Z_i := \{(x_1, \ldots, x_r) \ | \ x_i = x_{p(i)} \}.$$
We denote the class of $Z_i$ in the divisor class group by $\xi_i$ for any $i \in [1,r]$. Due to \cite[Proposition 3.5]{LT04} they form a basis of the cone of effective divisors. They even generate the Chow ring due to \cite{Dem74}, but we will not need this fact.

A basis $\mathcal{L}_1, \cdots, \mathcal{L}_r$ of the nef cone of the Bott-Samelson variety was described in \cite[Section 2.3]{Per07}. They are defined as follows. For any $i \in [1,r]$ there is a morphism $p_i: \widetilde{X}(\widetilde{w}) \to G/P^{\beta_i}$ induced by the projection. As $P^{\beta_i}$ is a maximal parabolic subgroup, the Picard group of $G/P^{\beta_i}$ is generated by a very ample line bundle $\mathcal{O}(1)$. We define $\mathcal{L}_i$ to be the pullback $p_i^* \mathcal{O}(1)$.

Similarly, we define line bundles on the intermediate variety. There is a morphism $\widehat{p_i}: \widehat{X}(\widehat{w}) \to G/P^{\beta_i}$ induced by the corresponding projection for all $i \in m_{\widehat{w}}(Q_w)$. We define the sheaf $\mathcal{L}^{\widehat{w}}_i$ as $\widehat{p}_i^*{\mathcal{O}(1)}$. For all $i \in p(Q_w)$ we denote $\widehat{D}_i$ as the pushforward $\widetilde{\pi}_*(\xi_i)$. Recall that the divisor class group of a Schubert variety has a basis given by the classes of Schubert divisors, which is also a basis of the cone of effective divisors. The following proposition is a consequence of the structure of the intermediate variety as a tower of locally trivial fibrations (see \cite{Per07} for more details).

\begin{prop}
\label{divInterm}
The sheaves $\mathcal{L}^{\widehat{w}}_i$ for $i \in m_{\widehat{w}}(Q_w)$ form a basis of the Picard group that generates the nef cone. Moreover, the divisors $\widehat{D}_i$ for all peaks $i \in p(Q_w)$ form a basis of the divisor class group and the cone of effective divisors. In particular, $\widehat{X}(\widehat{w})$ is $\mathbb{Q}$-factorial.
\end{prop}

The variety $\widehat{X}(\widehat{w})$ is $\mathbb{Q}$-factorial if and only if
the dimension of the space of $\mathbb{Q}$-divisors equals the rank of the
Picard group. The proposition says that this happens if and only if the number
of peaks is equal to the number of minimal vertices in the decomposition
$\widehat{w}$. This is true for all the decompositions described in this chapter
because every subquiver $Q_{w_i}$ has a unique peak.

\begin{lem}
\label{Li}
Let $\lambda^k_i$ be defined by
$$\mathcal{L}_i = \sum_{k=1}^r \lambda^k_i \xi_k.$$
Then we have $\lambda^k_i \geq 0$ for all $i,k$ and $\lambda^k_i = 0$ unless
$k \succeq i$.
\begin{proof}
The inequality $\lambda^k_i \geq 0$ follows directly from the definition of
$\mathcal{L}_i$ as the pullback of a very ample line bundle. The equality
$\lambda^k_i = 0$ unless $k \succeq i$ has been proven in \cite[Proposition
2.16]{Per07}.
\end{proof}
\end{lem}

\section{$\mathbb{Q}$-Factorialization}
\label{mori}

In this section we describe all $\mathbb{Q}$-factorializations of (co)minuscule
Schubert varieties in any characteristic. The proof relies on the theory of Mori
dream spaces. For more information on this topic we refer to \cite{HK00}.
Recall, that a birational morphism of varieties $f: Y \to X$ is called
Mori-small if it contracts no divisors.

\begin{defi}
Let $X$ be a normal projective variety. A $\mathbb{Q}$-factorialization of $X$
is a Mori-small birational morphism $f: Y \to X$ such that $Y$ is projective,
normal and $\mathbb{Q}$-factorial variety.
\end{defi}

\begin{thm}
\label{qfact}
Let $X(w)$ be a (co)minuscule Schubert variety. Then all $\mathbb{Q}$-factorializations of $X(w)$ are given by the morphisms $\widehat{\pi}: \widehat{X}(\widehat{w}) \to X(w)$ obtained from an ordering of the peaks.
\end{thm}

The remainder of this chapter is devoted to proving the theorem. Let $\widehat{w} = (w_1, \ldots, w_s)$ be a generalized decomposition of an element $w \in W$ as described in Section \ref{prelim}. By $p(Q_w)$ we denote the set of peaks of $Q_w$. The images of the divisors $(\widehat{D}_i)_{i \in p(Q_w)}$ under the morphism $\widehat{\pi}$ are exactly the Schubert divisors in $X(w)$. Therefore, the map $\widehat{\pi}$ is indeed Mori-small. Recall that an intermediate variety is $\mathbb{Q}$-factorial due to Proposition \ref{divInterm}. This shows that the intermediate varieties are $\mathbb{Q}$-factorializations of the corresponding Schubert variety.

In order to get the converse statement, we look at Mori dream spaces that were introduced in \cite{HK00}. Let $X$ be a projective, normal and $\mathbb{Q}$-factorial variety. A \textbf{small $\mathbb{Q}$-factorial modification (SQM)} of $X$ is a birational map $f: X \dashrightarrow Y$ inducing an isomorphism in codimension $1$ such that $Y$ is again a projective, normal and $\mathbb{Q}$-factorial variety.

\begin{defi}
\label{moriDream}
Let $X$ be a projective, normal and $\mathbb{Q}$-factorial variety. Then $X$ is called a \textbf{Mori dream space} if there is a finite set of SQMs $f_i: X \dashrightarrow X_i$ such that
\begin{enumerate}
\item the Picard group of $X$ is finitely generated,
\item the cone $\operatorname{Nef}(X_i)$ is polyhedral and generated by finitely many semiample divisors for all $i$,
\item the cone of movable divisors $\operatorname{Mov}(X)$ is the union of the $f_i^*(\operatorname{Nef}(X_i))$.
\end{enumerate}
\end{defi}

These spaces are tailor-made for running the minimal model program regardless of the characteristic of the field. We will only use the following statement from \cite[Proposition 1.11]{HK00}.

\begin{prop}
The SQMs used in the definition of a Mori dream space $X$ are all possible SQMs of $X$.
\end{prop}

\begin{prop}
Let $\widehat{w}=\widehat{w}^1, \widehat{w}^2, \ldots, \widehat{w}^t$ be all
possible generalized reduced decompositions of $w$ obtained from an ordering of
the peaks $p(Q_w)$.
\begin{enumerate}
  \item The birational maps $f_j: \widehat{X}(\widehat{w}) \dashrightarrow
  \widehat{X}(\widehat{w}^j)$ are small $\mathbb{Q}$-factorial modifications.
  \item The variety $\widehat{X}(\widehat{w})$ is a Mori dream
  space with small $\mathbb{Q}$-factorial modifications $f_j: \widehat{X}(\widehat{w}) \dashrightarrow
  \widehat{X}(\widehat{w}^j)$.
\end{enumerate}

\begin{proof}
Part (i) follows from the description of divisors in Proposition
\ref{divInterm}. The fact that $\operatorname{Pic}(\widehat{X}(\widehat{w}))$ is
finitely generated and that the nef cone is polyhedral with finitely many
semiample generators is stated in the same proposition.

Therefore, for all peaks $i \in p(Q_w)$ we can identify the divisors
$\widehat{D}_i$ in $\widehat{X}(\widehat{w})$ and in
$\widehat{X}(\widehat{w}^j)$ via $f_j^*$. The proof can be concluded by showing
the equalities $$\operatorname{Eff}(\widehat{X}(\widehat{w})) =
\operatorname{Mov}(\widehat{X}(\widehat{w})) = \bigcup_j
\operatorname{Nef}(\widehat{X}(\widehat{w}^j)).$$ The cone of effective divisors
is polyhedral with basis given by $(\widehat{D}_i)_{i \in p(Q_w)}$. Therefore,
it is closed in this case which proves
$\operatorname{Mov}(\widehat{X}(\widehat{w})) \subset
\operatorname{Eff}(\widehat{X}(\widehat{w}))$ and
$\operatorname{Nef}(\widehat{X}(\widehat{w}^j)) \subset
\operatorname{Eff}(\widehat{X}(\widehat{w}))$ for all $j$.

Let $D$ be any effective divisor in $\widehat{X}(\widehat{w})$. Then there are
coefficients $\lambda_p \geq 0$ such that $D = \sum_{p \in p(Q_w)} \lambda_p
\widehat{D}_p$. We define $N:= \{ p \in p(Q_w) \ | \ \lambda_p = 0 \}$, $t := \#
p(Q_w) \backslash N$ and $s := \# p(Q_w)$. We will prove the following by
induction on $t$. For any ordering $(p_{t+1}, \ldots, p_s)$ of $N$ there is an
ordering $(p_1, \ldots, p_t)$ of $p(Q_w)\backslash N$ such that for
$\widehat{w}^j$ obtained from $(p_1, \ldots, p_s)$, the divisor $D$ is in the
cone spanned by $(\mathcal{L}_{i}^{\widehat{w}^j})_{i \in m_{\widehat{w}^j}(Q_w)}$.

The case $t = 0$ is obvious. For $t>0$ we define
$$Q':= \{i \in Q_w \ | \ \nexists p \in N: \ i \preceq p \}.$$
Let $i_0$ be a minimal vertex of $Q'$. The coefficients of
$\mathcal{L}_{i_0}^{\widehat{w}^j}$ with respect to the base $(\widehat{D}_p)_{p
\in p(Q_w)}$ are independent from an ordering as above because of Lemma \ref{Li} and
the fact that $\mathcal{L}_{i_0}^{\widehat{w}^j}$ is the pushforward of
$\mathcal{L}_{i_0}$.
The same lemma also implies the existence of a coefficient
$\mu \in \mathbb{Q}_{\geq 0}$ such that $$D':= D - \mu
\mathcal{L}_{i_0}^{\widehat{w}^j} = \sum_{p \in p(Q_w)\backslash N} \lambda'_p
\widehat{D}_p$$ is effective and there is a peak $p_t$ such that $\lambda'_{p_t}
= 0$. As the inequality $t>\# \{p \in p(Q_w) \ | \ \lambda'_p \neq 0 \}$ holds,
we can conclude by induction.

This implies the remaining two inclusions
$$\operatorname{Eff}(\widehat{X}(\widehat{w})) \subset \operatorname{Mov}(\widehat{X}(\widehat{w})) \text{ and } \operatorname{Eff}(\widehat{X}(\widehat{w})) \subset \bigcup_j \operatorname{Nef}(\widehat{X}(\widehat{w}^j)).$$

\end{proof}
\end{prop}

Theorem \ref{qfact} follows immediately because any $\mathbb{Q}$-factorialization of $X(w)$ provides an SQM to an intermediate variety $\widehat{X}(\widehat{w})$.

\section{IH-Small Resolutions}
\label{IHresolution}

In this chapter we show how to derive a classification of all IH-small resolutions of (co)minuscule Schubert varieties from the classification of $\mathbb{Q}$-factorializations. The minuscule case was proven over $\mathbb{C}$ in \cite{Per07}.

\begin{defi}
\label{IHsmall}
A projective birational morphism $\pi: Y \to X$ between normal varieties is called \textbf{IH-small} if for all $k>0$ the following inequality holds.
$$\operatorname{codim}_X \{ x \in X \ | \ \dim \pi^{-1}(x) = k \} > 2k$$
In addition, if $Y$ is smooth, then $\pi$ is called an IH-small resolution.
\end{defi}

It is easy to see that any IH-small resolution is a
$\mathbb{Q}$-factorialization. Therefore, we only need to figure out which of
the $\mathbb{Q}$-factorializations are smooth and IH-small. In order to check
smoothness, we can use the following result due to \cite{BP99}.

\begin{prop}
Let $X(w)$ be any (co)minuscule Schubert variety that is not minuscule in the $C_n$-case. Then $X(w)$ is smooth if and only if it is homogeneous under its stabilizer.
\end{prop}

Notice that the only minuscule Schubert varieties in the $C_n$-case are projective spaces which have a bigger general linear group acting on them. The stabilizer in the symplectic group might not be big enough to make the variety homogeneous under its action.

A straightforward proof shows that the stabilizer $P_w$ of $X(w)$ is generated by $B$ and the sets $Bs_{\alpha}B$ whenever $\overline{s_{\alpha} w} \leq \overline{w}$ in the Bruhat order of $W/W_P$. As the varieties $\widehat{X}(\widehat{w})$ are towers of locally trivial fibration with fibers being (co)minuscule Schubert varieties, this can be used to show that $\widehat{X}(\widehat{w})$ is smooth if and only if every part of the decomposition of $Q_w$ corresponds to a smooth Schubert variety.

\begin{defi}
Let $w \in W$ be (co)minuscule. A vertex $h \in Q_w$ is called a \textbf{hole} if it has no predecessor and $\beta_h w > w$, i.e., gluing $\beta_h$ on top of the quiver $Q_w$ still leads to a quiver satisfying the conditions of Theorem \ref{propQuiv2}.
\end{defi}

\begin{prop}
A (co)minuscule Schubert variety $X(w)$ is homogeneous under its stabilizer, i.e., smooth, if and only if the quiver $Q_w$ has no holes. 
\end{prop}
This was proven in \cite[Theorem 0.2]{Per09} using the result of Brion and Polo. Note that Perrin denotes our holes as essential holes.

\begin{ex}
The vertices $h_1$, \ldots, $h_6$ are all holes in the following three examples. In fact, none of the corresponding Schubert varieties are smooth.

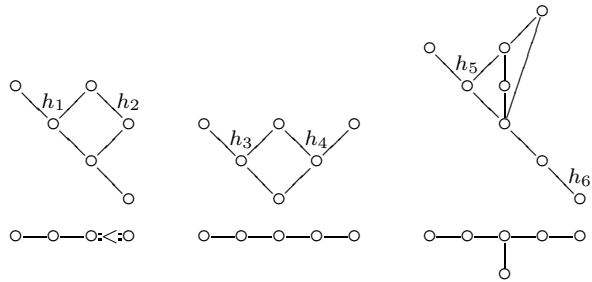
\centerline{
\xygraph{
!{<0cm,0cm>;<1cm,0cm>:<0cm,1cm>::}
!{(0,-1) }*{\circ}="1"
!{(1,-1) }*{\circ}="2"
!{(0.5,-1.5) }*{\circ}="3"
!{(1.5,-1.5) }*{\circ}="4"
!{(1,-2) }*{\circ}="5"
!{(1.5,-2.5) }*{\circ}="6"
!{(0,-3) }*{\circ}="a"
!{(0.5,-3) }*{\circ}="b"
!{(1,-3) }*{\circ}="c"
!{(1.5,-3) }*{\circ}="d"
!{(0.5,-1.2) }*{\scriptstyle{h_1}}="h1"
!{(1.5,-1.2) }*{\scriptstyle{h_2}}="h2"
!{(2.5,-1.5) }*{\circ}="1'"
!{(3.5,-1.5) }*{\circ}="2'"
!{(4.5,-1.5) }*{\circ}="3'"
!{(3,-2) }*{\circ}="4'"
!{(4,-2) }*{\circ}="5'"
!{(3.5,-2.5) }*{\circ}="6'"
!{(2.5,-3) }*{\circ}="a'"
!{(3,-3) }*{\circ}="b'"
!{(3.5,-3) }*{\circ}="c'"
!{(4,-3) }*{\circ}="d'"
!{(4.5,-3) }*{\circ}="e'"
!{(3,-1.7) }*{\scriptstyle{h_3}}="h3"
!{(4,-1.7) }*{\scriptstyle{h_4}}="h4"
!{(7,0) }*{\circ}="1''"
!{(5.5,-0.5) }*{\circ}="2''"
!{(6.5,-0.5) }*{\circ}="3''"
!{(6,-1) }*{\circ}="4''"
!{(6.5,-1) }*{\circ}="5''"
!{(6.5,-1.5) }*{\circ}="6''"
!{(7,-2) }*{\circ}="7''"
!{(7.5,-2.5) }*{\circ}="8''"
!{(5.5,-3) }*{\circ}="a''"
!{(6,-3) }*{\circ}="b''"
!{(6.5,-3) }*{\circ}="c''"
!{(7,-3) }*{\circ}="d''"
!{(7.5,-3) }*{\circ}="e''"
!{(6.5,-3.5) }*{\circ}="f''"
!{(6,-0.7) }*{\scriptstyle{h_5}}="h5"
!{(7.5,-2.2) }*{\scriptstyle{h_6}}="h6"
"a"-"b"
"b"-"c"
"c"-@{=}|{<}"d"
"1"-"3"
"2"-"3"
"2"-"4"
"3"-"5"
"4"-"5"
"5"-"6"
"a'"-"b'"
"b'"-"c'"
"c'"-"d'"
"d'"-"e'"
"1'"-"4'"
"2'"-"4'"
"2'"-"5'"
"3'"-"5'"
"4'"-"6'"
"5'"-"6'"
"a''"-"b''"
"b''"-"c''"
"c''"-"d''"
"c''"-"f''"
"d''"-"e''"
"1''"-"3''"
"1''"-"6''"
"2''"-"4''"
"3''"-"4''"
"3''"-"5''"
"4''"-"6''"
"5''"-"6''"
"6''"-"7''"
"7''"-"8''"
}}
\end{ex}

Recall the definition of the height $h(i)$ of a vertex $i$ of the quiver $Q_w$
as the largest integer $n$ such that there is a path consisting of $n-1$ arrows
from $i$ to the unique maximal vertex in $Q_w$. An ordering of the peaks $(p_1,
\ldots, p_{s+1})$ is called \textbf{neat} if $h(p_i) \leq h(p_{i+1})$ for all $i
\in [1,s]$. We should point out that this definition is slightly more
restrictive than the one given by Perrin in \cite{Per07} and by Zelevinsky in
\cite{Zel83}. However, it is not difficult to see that up to commutation in the
resulting generalized decompositions, every neat ordering of Perrin can be
modified to be of the above form. In particular, the actual intermediate
varieties are not affected. Perrin needed the slightly more complicated
definition for notational reasons in his proofs. We prefer this much easier
definition for our purposes.

\begin{cor}
\label{mainThm}
Let $X(w)$ be a (co)minuscule Schubert variety over an algebraically closed field $k$. Then the IH-small resolutions of $X(w)$ are exactly given by the morphisms $\widehat{\pi}: \widehat{X}(\widehat{w}) \to X(w)$, where $\widehat{w}$ is obtained from a neat ordering of the peaks and $\widehat{X}(\widehat{w})$ is smooth.
\begin{proof}
The fact that all resolutions coming from such an ordering are IH-small can be proven as in \cite{SV94} and in \cite{Per07}. We will not repeat the proof. Every IH-small morphism is a $\mathbb{Q}$-factorialization. By Theorem \ref{qfact} we are left to show that the other smooth varieties $\widehat{X}(\widehat{w})$ do not yield an IH-small resolution. This is done in an explicit construction in the proof of Proposition 7.2 of \cite{Per07}.
\end{proof}
\end{cor}

\begin{ex}
In our standard examples the first decomposition yields a non smooth variety. The second decomposition corresponds to an IH-small resolution. The third decomposition corresponds to a smooth variety, but the decomposition is not neat.

\centerline{
\xygraph{
!{<0cm,0cm>;<1cm,0cm>:<0cm,1cm>::}
!{(0,-1) }*{\circ}="1"
!{(1,-1) }*{\circ}="2"
!{(0.5,-1.5) }*{\circ}="3"
!{(1.5,-1.5) }*{\circ}="4"
!{(1,-2) }*{\circ}="5"
!{(1.5,-2.5) }*{\circ}="6"
!{(0,-3) }*{\circ}="a"
!{(0.5,-3) }*{\circ}="b"
!{(1,-3) }*{\circ}="c"
!{(1.5,-3) }*{\circ}="d"
!{(0.5,-1) }*{}="cut1"
!{(1.5,-2) }*{}="cut2"
!{(2.5,-1.5) }*{\circ}="1'"
!{(3.5,-1.5) }*{\circ}="2'"
!{(4.5,-1.5) }*{\circ}="3'"
!{(3,-2) }*{\circ}="4'"
!{(4,-2) }*{\circ}="5'"
!{(3.5,-2.5) }*{\circ}="6'"
!{(2.5,-3) }*{\circ}="a'"
!{(3,-3) }*{\circ}="b'"
!{(3.5,-3) }*{\circ}="c'"
!{(4,-3) }*{\circ}="d'"
!{(4.5,-3) }*{\circ}="e'"
!{(3,-1.5) }*{}="cut1'"
!{(3.5,-2) }*{}="cut2'"
!{(4,-1.5) }*{}="cut3'"
!{(3,-2.5) }*{}="cut4'"
!{(7,0) }*{\circ}="1''"
!{(5.5,-0.5) }*{\circ}="2''"
!{(6.5,-0.5) }*{\circ}="3''"
!{(6,-1) }*{\circ}="4''"
!{(6.5,-1) }*{\circ}="5''"
!{(6.5,-1.5) }*{\circ}="6''"
!{(7,-2) }*{\circ}="7''"
!{(7.5,-2.5) }*{\circ}="8''"
!{(5.5,-3) }*{\circ}="a''"
!{(6,-3) }*{\circ}="b''"
!{(6.5,-3) }*{\circ}="c''"
!{(7,-3) }*{\circ}="d''"
!{(7.5,-3) }*{\circ}="e''"
!{(6.5,-3.5) }*{\circ}="f''"
!{(5.75,-0.5) }*{}="cut1''"
!{(6.75,-1.5) }*{}="cut2''"
"cut1"-@{.}"cut2"
"a"-"b"
"b"-"c"
"c"-@{=}|{<}"d"
"1"-"3"
"2"-"3"
"2"-"4"
"3"-"5"
"4"-"5"
"5"-"6"
"cut1'"-@{.}"cut2'"
"cut3'"-@{.}"cut4'"
"a'"-"b'"
"b'"-"c'"
"c'"-"d'"
"d'"-"e'"
"1'"-"4'"
"2'"-"4'"
"2'"-"5'"
"3'"-"5'"
"4'"-"6'"
"5'"-"6'"
"cut1''"-@{.}"cut2''"
"a''"-"b''"
"b''"-"c''"
"c''"-"d''"
"c''"-"f''"
"d''"-"e''"
"1''"-"3''"
"1''"-"6''"
"2''"-"4''"
"3''"-"4''"
"3''"-"5''"
"4''"-"6''"
"5''"-"6''"
"6''"-"7''"
"7''"-"8''"
}}
\end{ex}

\renewcommand{\refname}{References}

{\sc Department of Mathematics, The Ohio State University, 231 W 18th Avenue, Columbus, OH 43210-1174, USA}

{\it E-mail address:} {\tt schmidt.707@osu.edu}

{\it URL:} {\tt https://people.math.osu.edu/schmidt.707/}

\end{document}